\documentclass[11pt]{amsart}

\usepackage{latexsym}
\usepackage{amscd}
\usepackage[all]{xy}
\usepackage[mathscr]{euscript}

\normalsize

\RequirePackage{xspace}

\newcommand{\thought}[1]{}
\renewcommand{\thought}[1]{ \textbf{[#1]}}

\usepackage{enumerate}        
\newenvironment{roenumerate}{\begin{enumerate}[\upshape (i)]}{\end{enumerate}}

\newcommand\nc {\newcommand}
\newcommand\rnc{\renewcommand}

\newcount\blopone
\newcount\xone
\newcount\xtwo
\newcount\ytwo

\newtheorem{theorem}{Theorem}[section]
\newtheorem{prop}[theorem]{Proposition}
\newtheorem{com}[theorem]{Comment}
\newtheorem{redu}[theorem]{Reduction}
\newtheorem{refinement}[theorem]{Refinement}
\newtheorem{summary}[theorem]{Summary}
\newtheorem{importnota}[theorem]{Important Notation}
\newtheorem{prblm}[theorem]{Problem}
\newtheorem{notation}[theorem]{Notation}
\newtheorem{defin}[theorem]{Definition}
\newtheorem{caution}[theorem]{Caution}
\newtheorem{remark}[theorem]{Remark}
\newtheorem{reminder}[theorem]{Reminder}
\newtheorem{illustration}[theorem]{Illustration}
\newtheorem{lemma}[theorem]{Lemma}
\newtheorem{construction}[theorem]{Construction}
\newtheorem{corollary}[theorem]{Corollary}
\newtheorem{example}[theorem]{Example}
\newtheorem{conclusion}[theorem]{Conclusion}
\newtheorem{triviality}[theorem]{Triviality}
\newtheorem{proto}[theorem]{Prototype Quasifibration}
\newtheorem{cauex}[theorem]{Cautionary Example}
\newtheorem{hypo}[theorem]{Hypothesis}
\newtheorem{subth}{ }[theorem]
\newtheorem{case}{Case}[theorem]
\newtheorem{ssubth}{ }[subth]
\newtheorem{facts}[theorem]{Facts}

\nc\tri[1]{\begin{triviality}
\label{#1}}
\nc\fac[1]{\begin{facts}
\label{#1}
\begin{em}}
\nc\cas[1]{\begin{case}
\label{#1}
\begin{em}}
\nc\rfn[1]{\begin{refinement}
\label{#1}}
\nc\prt[1]{\begin{proto}
\label{#1}}
\nc\lem[1]{\begin{lemma}
\label{#1}}
\nc\pro[1]{\begin{prop}
\label{#1}}
\nc\thm[1]{\begin{theorem}
\label{#1}}
\nc\cor[1]{\begin{corollary}
\label{#1}}
\nc\dfn[1]{\begin{defin}
\label{#1}}
\nc\sthm[1]{\begin{subth}
\label{#1}}
\nc\exm[1]{\begin{example}
\label{#1}
\begin{em}}
\nc\plm[1]{\begin{prblm}
\label{#1}
\begin{em}}
\nc\rmk[1]{\begin{remark}
\label{#1}
\begin{em}}
\nc\rmd[1]{\begin{reminder}
\label{#1}
\begin{em}}
\nc\ntn[1]{\begin{notation}
\label{#1}
\begin{em}}
\nc\smr[1]{\begin{summary}
\label{#1}
\begin{em}}
\nc\cau[1]{\begin{caution}
\label{#1}
\begin{em}}
\nc\hyp[1]{\begin{hypo}
\label{#1}}
\nc\imn[1]{\begin{importnota}
\label{#1}
\begin{em}}
\nc\rdn[1]{\begin{redu}
\label{#1}
\begin{em}}
\nc\cax[1]{\begin{cauex}
\label{#1}
\begin{em}}
\nc\cmt[1]{\begin{com}
\label{#1}
\begin{em}}
\nc\con[1]{\begin{construction}
\label{#1}
\begin{em}}
\nc\ill[1]{\begin{illustration}
\label{#1}
\begin{em}}
\nc\ssthm[1]{\begin{ssubth}
\label{#1}
\begin{em}}
\nc\cnc[1]{\begin{conclusion}
\label{#1}
\begin{em}}

\nc\elem{\end{lemma}}
\nc\erdn{\end{em}\end{redu}}
\nc\erfn{\end{refinement}}
\nc\eprt{\end{proto}}
\nc\ethm{\end{theorem}}
\nc\ecor{\end{corollary}}
\nc\edfn{\end{defin}}
\nc\esthm{\end{subth}}
\nc\epro{\end{prop}}
\nc\etri{\end{triviality}}
\nc\eexm{\end{em}
\end{example}}
\nc\ecmt{\end{em}
\end{com}}
\nc\efac{\end{em}
\end{facts}}
\nc\ermk{\end{em}
\end{remark}}
\nc\ermd{\end{em}
\end{reminder}}
\nc\eill{\end{em}
\end{illustration}}
\nc\eplm{\end{em}
\end{prblm}}
\nc\ecas{\end{em}
\end{case}}
\nc\ecau{\end{em}
\end{caution}}
\nc\ecax{\end{em}
\end{cauex}}
\nc\eimn{\end{em}
\end{importnota}}
\nc\entn{\end{em}
\end{notation}}
\nc\econ{\end{em}
\end{construction}}
\nc\esmr{\end{em}
\end{summary}}
\nc\ehyp{
\end{hypo}}
\nc\ecnc{\end{em}
\end{conclusion}}
\nc\essthm{\end{em}
\end{ssubth}}

\nc\sst{\scriptstyle}
\newcommand{\comment}[1]{}
\newcommand{\ri}{\longrightarrow}

\newcommand{\K}{{\mathbf K}}
\newcommand{\D}{{\mathbf D}}

\nc\op{^{\hbox{\rm\tiny op}}}
\nc\mth{^{\hbox{\rm\tiny th}}}

\nc\script{\mathscr}
\nc\z{\zeta}
\nc\bc{{\mathbb{BC}}}
\nc\ct{{\script T}}
\nc\cf{{\script F}}
\nc\cl{{\script L}}
\nc\cv{{\script V}}
\nc\ce{{\script E}}
\nc\cs{{\script S}}
\nc\car{{\script R}}
\nc\cd{{\script D}}
\nc\cc{{\script C}}
\nc\ca{{\script A}}
\nc\ci{{\script I}}
\nc\co{{\script O}}
\nc\bd{\begin{description}}
\nc\ed{\end{description}}
\nc\ctob{{\script C}at\big(\ci^{op},\ca\big)}
\nc\clim{{\ds\mathop{\rm lim}_{\ds\longleftarrow}}\,}
\nc\climi{\clim^{\!i}\,}
\nc\climn{\clim^{\!n}\,}
\nc\colim{{\ds\mathop{\rm colim}_{\ds\la}}}
\nc\oa{\overline{\ca}}
\nc\s{\sigma}
\nc\ta{\tau}
\nc\os{\overline\sigma}
\nc\ot{\overline\tau}
\nc\T{\Sigma}
\nc\de[1]{{\mathop{\rm deg(#1)}}}
\nc\Ad[1]{\mathop{\rm Ad}(#1)}
\nc\ad[1]{\mathop{\rm ad}(#1)}
\nc\kth{{\it K}--theory}
\nc\loc[1]{{\text{\rm Loc(#1)}}}
\nc\coloc[1]{{\text{\rm Coloc}(#1)}}

\def\der #1 {D\left(#1\right)}
\nc\prf{\begin{proof}}
\nc\eprf{\end{proof}}
\nc\ds{\displaystyle}
\nc\Tor{\text{\rm Tor}}

\nc\cb{{\script B}}
\nc\ab{{\script A}b}

\nc\be{\begin{roenumerate}}
\nc\ee{\end{roenumerate}}

\nc\cat[1]{{\script C}at\Big({\big\{#1\big\}}\op\,\,,\,\,\ab\Big)}
\nc\csab{{\script C}at\big(\cs^{op},\ab\big)}
\nc\ctab{{\script C}at\Big({\{\ct^\alpha\}}^{op},\ab\Big)}
\nc\csex{{\script E}x\big(\cs^{op},\ab\big)}
\nc\ctex{{\script E}x\Big({\{\ct^\alpha\}}^{op},\ab\Big)}
\nc\sub{\qquad\subset\qquad}
\nc\ctr[1]{{\left.\ct\left(-,#1\right)\right|}_{\cs}}
\nc\ctrf[2]{{\left.\ct\left(#1,#2\right)\right|}_{\cs}}
\nc\Ctr[1]{{\left.\ct\left(-,#1\right)\right|}_{\ct^\alpha}}
\nc\Ctrf[2]{{\left.\ct\left(#1,#2\right)\right|}_{\ct^\alpha}}

\nc\la{\longrightarrow}
\nc\nin{\noindent}
\nc\cad[1]{\text{card}(#1)}
\nc\eq{\quad=\quad}
\nc\BA{\begin{array}{c}}
\nc\EA{\end{array}}
\nc\barr{
\[
\begin{array}{cccccccccccccccc}
}
\nc\earr{
\end{array}
\]
}
\nc\as[1]{{\langle S\rangle}^{#1}}
\nc\sh{\text{\it shift}}

\nc\yy[1]{{\left.\ct\left(-,#1\right)\right|}_{\ct^c}}
\nc\vrep[2]{{\left.\ct\left(#1,#2\right)\right|}_{\ct^\alpha}}
\nc\da{\downarrow}
\nc\Hom{{\mathop{\rm Hom}}}
\nc\HHom{{\script H}{\mathop{\rm om}}}
\nc\End{{\mathop{\rm End}}}
\nc\Ext{{\mathop{\rm Ext}}}
\nc\mr\Modtc
\nc\PExt{{\mathop{\rm PExt}}}
\nc\stm{\text{\rm stmod}(kG)}
\nc\stM{\text{\rm StMod}(kG)}
\nc\e{\varepsilon}
\nc\p{\varphi}

\nc\rs{\s^{-1}A}
\nc\br{{\{\s^{-1}A\}}}
\nc\ra\ri
\nc\y[1]{\mathbf{y}#1}
\nc\x[1]{\mathbf{z}#1}
\nc\mmod[1]{#1\text{--\rm mod}}
\nc\Mod[1]{#1\text{--\rm Mod}}
\nc\Md {\ensuremath{\mathop{\textup{Mod}}}}
\rnc\mod[1]{\ensuremath{\mathop{\textup{mod-}#1}}\xspace}
\nc\Modtc{\Mod{\ct^c}}
\nc\pgldim[1]{\mathop{\rm pgldim}\,#1}
\nc\tf{{\rm [TR5]}}
\nc\tfs{{\rm [TR5$^*$]}}
\nc\Fun{\text{\rm Funct}(F\op,\ab)}
\nc\sym{\text{\rm Sym}}
\nc\sgn{\text{\rm sgn}}
\nc\Pro{\text{\rm Prod}^{}_\alpha(F\op,\ab)}
\nc\Yt[1]{{\left.\Hom_\ct^{}\left(-,#1\right)\right|}_F^{}}
\nc\dl{\delta}
\nc\Proj[1]{#1\text{--\rm Proj}}
\nc\proj[1]{#1\text{--\rm proj}}
\nc\Flat[1]{#1\text{--\rm Flat}}
\nc\Inj[1]{#1\text{--\rm Inj}}
\nc\ov{\overline}
\nc\wt{\widetilde}
\nc\ph{\varphi}
\nc\tstr{{\it t}--structure}
\nc\spec[1]{{\text{\rm Spec}(#1)}}
\newcommand{\m}{\mathfrak{m}}
\newcommand{\n}{\mathfrak{n}}

\nc\hoco{
\begin{picture}(40,10)
\put(20,0){\makebox(0,0)[b]{\text{\rm Hocolim}}}
\put(5,-2){\vector(1,0){30}}
\end{picture}\,\,}

\nc\holim{
\begin{picture}(40,10)
\put(20,0){\makebox(0,0)[b]{\text{\rm Holim}}}
\put(35,-2){\vector(-1,0){30}}
\end{picture}}

\begin{document}

\author{Amnon Neeman}\thanks{The research was partly supported 
by the Australian Research Council}
\address{Centre for Mathematics and its Applications \\
        Mathematical Sciences Institute\\
        John Dedman Building\\
        The Australian National University\\
        Canberra, ACT 0200\\
        AUSTRALIA}
\email{Amnon.Neeman@anu.edu.au}

\title[Cogenerators in $\K(\Proj R)$]{Non-left-complete derived categories}

\begin{abstract}
We give some examples of abelian categories $\ca$ for which the
derived category $\D(\ca)$ is not left-complete. Perhaps the most
natural of these is where $\ca$ is the category of representations
of the additive group $\mathbb{G}_a$ over a field $k$ of
characteristic $p>0$.
\end{abstract}

\subjclass[2000]{Primary 18E30, secondary 18G55}

\keywords{Derived categories, {\it t}--structures, homotopy limits}

\maketitle

\tableofcontents

\setcounter{section}{-1}

\section{Assumed background}
\label{S0}

In this article we assume the reader is familiar with derived categories 
and with {\it t}--structures on them. See Verdier~\cite{Verdier96} for 
the theory of derived categories, and Beilinson, Bernstein and 
Deligne~\cite[Chapter~1]{BeiBerDel82} for an introduction to 
{\it t}--structures.

\section{The counterexample}
\label{S1}

Suppose $\ca$ is an abelian category and $\D(\ca)$ is its derived category.
For any object $x\in\D(\ca)$, we write $x^{\geq n}$ for the truncation
of $x$ with respect to the standard \tstr. We have canonical maps
$x^{\geq n}\la x^{\geq n+1}$, and a (non-canonical) map 
\[
\CD
\ph_x\,\,:\,\,x @>>> \holim x^{\geq n}\ .
\endCD
\]
The category $\D(\ca)$ is said to be \emph{left-complete} if, for every 
object $x\in\D(\ca)$, any map $\ph_x$ as above is an isomorphism. Even 
though the map $\ph_x$ is not canonical, it can be shown that, for given 
$x$, if one $\ph_x$ is an isomorphism then they all are.

The reader can find much more about left-complete categories in 
Lurie~\cite[Section~7]{Lurie06} 
or~\cite[Subsection~1.2.1, more precisely starting from 
Proposition~1.2.1.17]{Lurie11}. See also 
Drinfeld and Gaitsgory~\cite{DrinfeldGaitsgory11}.

In this note we will see how to produce many $\ca$
for which $\D(\ca)$ is not left-complete. Our counterexamples
will be of a very special form, which allows us to easily compute the
homotopy inverse limit $\holim x^{\geq n}$. Let us now sketch what we will do.

We will suppose that the 
abelian category $\ca$ satisfies the axiom [AB4], that is
coproducts are exact; this makes it easy to compute
coproducts in the derived category $\D(\ca)$, just form
the coproduct as complexes. Suppose $A$ is an object in our [AB4] 
abelian category $\ca$, and let 
\[
x\eq \coprod_{i=0}^\infty A[i]\ .
\]
It is clear that, for $n>0$, we have
\[
x^{\geq-n}\eq \coprod_{i=0}^n A[i]\eq \prod_{i=0}^n A[i]\ ,
\]
where the last equality is because finite coproducts agree with finite 
products. Now the homotopy inverse limit of the products is a genuine
inverse limit, and we have
\[
\holim x^{\geq n}\eq  \prod_{i=0}^\infty A[i]\ .
\]
Thus our problem becomes to decide whether the map
\[
\CD
\ds\coprod_{i=0}^\infty A[i]@>\ph>>   \ds\prod_{i=0}^\infty A[i]
\endCD
\]
is an isomorphism. Note that in this case the map is canonical; our 
homotopy inverse limit happens to be a genuine inverse limit, removing
the arbitrariness.
The left hand side is easy to work with; its cohomology is $A$ in each degree
$n\leq0$. What we will show is how to produce examples where the right hand
side has lots more cohomology. More precisely, we have
\[
\prod_{i=0}^\infty A[i]\eq A[0]\oplus\left(
\prod_{i=1}^\infty A[i]
\right)
\]
and the expectation would be for the second term to have a vanishing $H^0$;
what we will show is how to
produce non-zero 
classes in
\[
H^0\left(
\prod_{i=1}^\infty A[i]
\right)\ .
\]

It is time to disclose what will be our choice for the category $\ca$ and
for the object $A\in\ca$. 

\con{C0.1}
Let $k$ be a field, let $R_1$ be a finitely generated $k$ algebra, and 
let $\m$ be a $k$--point of $\spec {R_1}$. In other words, $\m\subset R$ 
is a maximal ideal with $R_1/\m\cong k$.
We make a string of definitions:
\be
\item
$
R_n=\otimes_{i=1}^n R_1
$, where the tensor is over the field $k$.
\item
The inclusion $R_n\la R_{n+1}$ is the inclusion of the tensor product
of the first $n$ terms.
\item
$R=\colim R_n$.
\item
The map $\Phi_i:R_1\la R$ is the inclusion of the $i\mth$ factor.
\item
The category $\ca$ will be the category of all those
$R$--modules, on which
$\Phi_i(\m)$ acts trivially for all but finitely many $i$.
\setcounter{enumiv}{\value{enumi}}
\ee
The object $A\in\ca$ will be the colimit over 
$n$ of the $R_n$--modules $k= \otimes_{i=1}^n [R_1/\m]$.
\econ

The main result is

\thm{T0.2}
Assume that $k=R_1/\m$ is not projective over the localization
${(R_1)}_\m$ of the ring $R_1$ at the maximal ideal $\m$.
With the category $\ca$ and the object $A\in\ca$ as in Construction~\ref{C0.1},
there is a non-zero element in
\[
H^0\left(
\prod_{i=1}^\infty A
\right)\ .
\]
\ethm

\rmk{R0.3}
The case where $R_1=k[x]/(x^p)$ is of particular interest. If the field 
$k$ is of characteristic $p$ then the category $\ca$ happens to be the category
of representations of the additive group $\mathbb{G}_a$, and we learn
that its derived category is not left-complete.
\ermk 

\rmk{boundunstable}
We trivially have
\[
\prod_{i=1}^\infty A[i]\eq \left(\prod_{i=1}^n A[i]\right)\oplus\left(
\prod_{i=n+1}^\infty A[i]\right)\ ,
\]
and hence
\[
H^0\left(
\prod_{i=1}^\infty A[i]
\right)\eq
H^0\left(
\prod_{i=1}^n A[i]
\right)\oplus
H^0\left(
\prod_{i=n+1}^\infty A[i]
\right)\ .
\]
On the other hand, with the finite product we have no problem computing
\[
H^0\left(
\prod_{i=1}^n A[i]
\right)\eq
H^0\left(
\coprod_{i=1}^n A[i]
\right)\eq0\ ,
\]
and Theorem~\ref{T0.2} now allows us to deduce that
\[
H^0\left(
\prod_{i=n+1}^\infty A[i]
\right)\quad\neq\quad
0\ .
\]
Translating we have that 
\[
H^n\left(
\prod_{i=1}^\infty A[i]
\right)\quad\neq\quad
0
\]
for all $n\geq 0$. The complexes $A[i],\, i>0$ all belong to 
${\D(\ca)}^{<0}$, but the
product $\prod_{i=1}^\infty A[i]$ is not bounded above.
\ermk

\medskip

\nin
{\bf Acknowledgements.}\ \ 
The author would like to thank Drinfeld and Gaitsgory for asking the question
that led to these counterexamples.

\section{The proof}
\label{S2}

We begin with a little lemma.

\lem{L2.1}
Let $k$ be a field, and let $R$ and $S$ be finitely generated $k$--algebras.
Suppose further that we are given $k$--points of $\spec R$ and $\spec S$; that
is $\m\subset R$ and $\n\subset S$ are maximal ideals, with
\[
R/\m\quad\cong\quad k\quad\cong\quad S/\n\ .
\]
Let $E$ be an injective envelope of $k=R/\m$ over the ring $R$, and 
$F$ an injective envelope of $k=S/\n$ over the ring $S$. Then $E\otimes_k F$
is an injective envelope of $k$ over the ring $R\otimes_k S$.
\elem

\prf
We will first prove the case where $R$ and $S$ are polynomial rings.

Let $R'=k[x_1^{},x_2^{},\ldots,x_m^{}]$ be a polynomial ring, and
let $\m$ be the maximal ideal generated by 
$\{x_1^{},x_2^{},\ldots,x_m^{}\}$. Then we know the injective envelope $E'$
of $k=R'/\m$ explicitly: it is the quotient of
$S=k[x_1^{},x_1^{-1},x_2^{},x_2^{-1},\ldots,x_m^{},x_m^{-1}]$ by the $R'$--submodule
generated by all monomials $x_1^{i_1^{}}x_2^{i_2^{}}\cdots x_m^{i_m^{}}$ with at 
least one
of the $i_j^{}>0$. As a $k$--vector space 
$E'=k[x_1^{-1},x_2^{-1},\ldots,x_m^{-1}]$,
and the $R'$--module structure is obvious when we declare 
$x_1^{i_1^{}}x_2^{i_2^{}}\cdots x_m^{i_m^{}}=0$ if 
some $i_j^{}>0$. If $S'=k[y_1^{},y_2^{},\ldots,y_n^{}]$ and $\n\subset S'$
is the ideal generated by $\{y_1^{},y_2^{},\ldots,y_n^{}\}$, then the fact that
\[
E'\otimes_k F'\eq k[x_1^{-1},x_2^{-1},\ldots,x_m^{-1}]\otimes_k
k[y_1^{-1},y_2^{-1},\ldots,y_n^{-1}]
\]
is the injective hull of $k$ over $R'\otimes S'$ is by inspection.

Now for the general case: assume $R=R'/I$ and $S=S'/J$ where $R'$ and $S'$ 
are polynomial rings, and $I\subset R'$ and $J\subset S'$ are ideals 
contained in the $\m$ and $\n$ above. Then the injective hull $E$ of $k=R/\m$
over the ring $R$ is the largest $R$--submodule of the $R'$--module
$E'$, that is
the $R'$--submodule $E\subset E'$ of all elements annihilated by the ideal $I$.
The lemma therefore comes down to the fact that the submodule of
$E'\otimes_k F'$ annihilated by the ideal 
$I\otimes_k S'+R'\otimes_k J$ is precisely
$E\otimes_k F$.
\eprf

\nin
{\bf Proof of Theorem~\ref{T0.2}.} \ \ 
Let $\ov R$ be the localization of $R_1$ at the maximal ideal $\m$.
We are assuming that $k$ is not projective over $\ov R$, that is the
projective dimension of $k$ is at least one.
Choose and fix a minimal free resolution of $k=\ov R/\m\ov R$ as an 
$\ov R$--module.
Let us write this resolution as 
\[
\CD
@>>>P_2@>>> P_1 @>>> P_0 @>>> k @>>> 0\ .
\endCD
\]
Then the modules $P_i$ are all finite and free over the ring $\ov R$, 
the differentials are
all matrices over $\ov R$, and the minimality guarantees that
the entries in these matrices all belong to the ideal $\ov\m=\m\ov R
\subset \ov R$.
Now let $E$ be the $\ov R$--injective envelope of the module $k$; applying the 
functor $\Hom_{\ov R}(-,E)$ to the projective resolution above, we produce
an injective resolution $I^*$ of $k$, which we write out as
\[
\CD
0@>>> k@>>>I^0@>>> I^1 @>>> I^2 @>>>
\endCD
\]
We know that each $I^j=\Hom(P_j,E)$ is a finite coproduct of copies of $E$,
and that the differentials $I^j\la I^{j+1}$ are matrices whose entries belong
to the ideal $\ov\m$. The fact that the projective dimension of $k$ is at least
one tells us that $P_1\neq 0$, and therefore $I^1\neq0$. Note
that an injective envelope $E$ of $k$ over the localized ring 
$\ov R={(R_1)}_\m$ is also an
injective envelope of $k$ over the ring $R_1$, hence we
have produced an injective resolution of $k$ over $R_1$. Next we
\be
\item
Choose a non-zero element $a$ in the image of the map $k\la I^0$.
\item
Choose a non-zero element $b\in I^1$, with $\m b=0$.
\setcounter{enumiv}{\value{enumi}}
\ee

If we view $k$ as a module over the ring $R_n=\otimes_{i=1}^nR_1$, then
the tensor product $J_n^*= \otimes_{i=1}^nI^*$ is certainly a resolution of
$k$ as an $R_n$ module, and Lemma~\ref{L2.1} guarantees further that
\be
\setcounter{enumi}{\value{enumiv}}
\item
Each $J_n^i$ is injective as a module over $R_n$.
\item
Let the inclusion $J_n^*\la J_{n+1}^*$ be the map
taking $x\in J_n^*$ to 
\[
x\otimes a\quad\in\quad
J_n^*\otimes I^0\quad\subset\quad
J_n^*\otimes I^*\eq J_{n+1}^*\ ,
\]
where $a\in I^0$ is as in (i) above. We define $J^*$ to be
\[
J^*\eq\colim J^*_n\ ;
\]
then $J^*$ is an injective resolution of $k$ in the category $\ca$.
\setcounter{enumiv}{\value{enumi}}
\ee

To prove the theorem we need to find a 
non-zero element in $H^0\left(\prod_{i>0} k[i]\right)$,
and our next observation is that the product
in the derived category $\prod_{i>0} k[i]$ is obtained as the
the ordinary product of injective
resolutions. The complex $J^*[i]$ is an injective
resolution of $k[i]$, and hence the derived product $\prod_{i>0} k[i]$ is
just the usual product 
$\prod_{i>0} J^*[i]$. Now for every $i\geq1$ let
\[S_i=\{i^2+1,\ldots,i^2+i\},\]
 and observe
that the sets $S_i$ are disjoint. In the injective $R_{i^2+i}$--module
\[
J^i_{i^2+i}=\coprod_{\sum \ell_m=i}I^{\ell_1}\otimes  I^{\ell_2}\otimes \cdots
\otimes I^{\ell_{i^2+i}}
\]
or more specifically in the summand 
\[
{(I^{0})}^{\otimes i^2}\otimes  {(I^{1})}^{\otimes i}
\]
we take the term 
\[
\lambda_i\eq a^{\otimes i^2}\otimes b^{\otimes i}\ ,
\] 
where $a\in I^0$ and 
$b\in I^1$ are as in (i) and (ii) above. The embedding $J^*_{i^2+i}\la J^*$
of (iv) gives us an element which we will denote $\lambda_i\in J^i$.
The elements
$\lambda_i$ have the properties
\be
\setcounter{enumi}{\value{enumiv}}
\item
Each $\lambda_i$ is a cycle; the differential $J^i\la J^{i+1}$ kills $\lambda_i$.
\item
$\Phi_j(\m)\lambda_i=0$
for all $i$ and $j$.
\setcounter{enumiv}{\value{enumi}}
\ee
We are assuming $i>0$, so 
each $\lambda_i$ must be a boundary because $H^i(J^*)=0$. But if
 $\mu_i\in J^{i-1}$ maps to $\lambda_i$, then there must exist an
 integer $j\in S_i$
so that $\Phi_j(\m)$ does not kill $\mu_i$.
Now form the element 
\[
\prod_{i=1}^\infty \lambda_i\quad\in\quad \prod_{i=1}^\infty J^i\ ,
\]
where the product is in the category of all $R$--modules.

\cau{C2.5}
The reader is reminded that the category $\ca$ is a subcategory of the 
category of $R$--modules. Both categories have infinite products; the products 
in the category of $R$--modules are just the usual cartesian products,
while the products in $\ca$ are subtler. To form the product in $\ca$
of a bunch of objects in $\ca$, one first forms the usual cartesian 
product, and then consider inside it the largest object belonging to $\ca$,
that is the collection of all elements satisfying part (v) of 
Construction~\ref{C0.1}.
\ecau

The element $\prod_{i=1}^\infty \lambda_i$
is a degree 0 cycle in the complex $\prod_{i\geq1}J^*[i]$, and it is annihilated
by $\Phi_j(\m)$ for all $j$. By Caution~\ref{C2.5} we have that
$\prod_{i=1}^\infty \lambda_i$ belongs to $\prod_{i=1}^\infty J^i$ even when
the product is understood in $\ca$. However, it is 
not a boundary in $\ca$.
If we try to express $\prod_{i=1}^\infty \lambda_i$ as the boundary
of
\[
\prod_{i=1}^\infty \mu_i\quad\in\quad \prod_{i=1}^\infty J^{i-1}\ ,
\]
then we discover that each $\mu_i$ fails to be annihilated by some $\Phi_j(\m)$
with $j\in S_i$. As the $S_i$ are disjoint, this produces infinitely
many $\Phi_j(\m)$ not annihilating $\prod_{i=1}^\infty \mu_i$, meaning
it does not belong to $\ca$.

\bibliographystyle{amsplain}
\bibliography{stan}

\def\cprime{$'$}
\providecommand{\bysame}{\leavevmode\hbox to3em{\hrulefill}\thinspace}
\providecommand{\MR}{\relax\ifhmode\unskip\space\fi MR }
\providecommand{\MRhref}[2]{%
  \href{http://www.ams.org/mathscinet-getitem?mr=#1}{#2}
}
\providecommand{\href}[2]{#2}
\begin{thebibliography}{1}

\bibitem{BeiBerDel82}
Alexander~A. Beilinson, Joseph Bernstein, and Pierre Deligne, \emph{Analyse et
  topologie sur les {\'e}spaces singuliers}, Ast{\'e}risque, vol. 100, Soc.
  Math. France, 1982 (French).

\bibitem{DrinfeldGaitsgory11}
Vladimir Drinfeld and Dennis Gaitsgory, \emph{On some finiteness questions for
  algebraic stacks}, in preparation.

\bibitem{Lurie06}
Jacob Lurie, \emph{Derived algebraic geometry {I}: stable $\infty$-categories},
  arXiv:math/0608228v5.

\bibitem{Lurie11}
\bysame, \emph{Higher {A}lgebra}, Prepint, available from
  http://www.math.harvard.edu/\hbox{$\sim$}lurie/.

\bibitem{Verdier96}
Jean-Louis Verdier, \emph{Des cat{\'e}gories d{\'e}riv{\'e}es des
  cat{\'e}gories abeliennes}, Asterisque, vol. 239, Soci{\'e}t{\'e}
  Math{\'e}matique de France, 1996 (French).

\end{thebibliography}

\end{document}